\input amstex

\input amstex
\input amsppt.sty
\magnification=\magstep1
\hsize=36truecc
 \vsize=23.9truecm
\baselineskip=14truept
 \NoBlackBoxes
\def\q{\quad}
\def\qq{\qquad}
\def\mod#1{\ (\text{\rm mod}\ #1)}

\def\t{\text}
\def\qtq#1{\q\t{#1}\q}
\def\mod#1{\ (\text{\rm mod}\ #1)}
\def\qtq#1{\q\t{#1}\q}
\def\f{\frac}
\def\e{\equiv}
\def\b{\binom}

\def\sls#1#2{(\f{#1}{#2})}
 \def\ls#1#2{\big(\f{#1}{#2}\big)}
\def\Ls#1#2{\Big(\f{#1}{#2}\Big)}
\par preprint: October 14, 2012
\par\q
\let \pro=\proclaim
\let \endpro=\endproclaim
\topmatter
\title Congruences concerning Legendre polynomials III
\endtitle
\author Zhi-Hong Sun\endauthor
\affil School of Mathematical Sciences, Huaiyin Normal University,
\\ Huaian, Jiangsu 223001, PR China
\\ Email: zhihongsun$\@$yahoo.com
\\ Homepage: http://www.hytc.edu.cn/xsjl/szh
\endaffil

 \nologo \NoRunningHeads

\abstract{Let $p>3$ be a prime, and let  $R_p$ be the set of
rational numbers whose denominator is coprime to $p$. Let
$\{P_n(x)\}$ be the Legendre polynomials. In this paper we mainly
show that for $m,n,t\in R_p$ with $m\not\e 0\pmod p$,
 $$\align &P_{[\frac
p6]}(t) \e -\Big(\frac
3p\Big)\sum_{x=0}^{p-1}\Big(\frac{x^3-3x+2t}p\Big)\pmod p,
\\&\Big(\sum_{x=0}^{p-1}\Big(\frac{x^3+mx+n}p\Big)\Big)^2\equiv
\Big(\frac{-3m}p\Big)
\sum_{k=0}^{[p/6]}\binom{2k}k\binom{3k}k\binom{6k}{3k}
\Big(\frac{4m^3+27n^2}{12^3\cdot 4m^3}\Big)^k\pmod p,\endalign
$$
 where $(\frac ap)$ is the Legendre symbol and $[x]$ is the
greatest integer function. As an application we solve some
conjectures of Z.W. Sun and the author concerning
$\sum_{k=0}^{p-1}\binom{2k}k\binom{3k}k\binom{6k}{3k}/m^k\pmod
{p^2}$, where $m$ is an integer not divisible by $p$.
\par\q
}
 \endabstract

\endtopmatter
\document
\subheading{1. Introduction}
\par Let $\{P_n(x)\}$ be the Legendre polynomials given by
$$P_0(x)=1,\ P_1(x)=x,\ (n+1)P_{n+1}(x)=(2n+1)xP_n(x)-nP_{n-1}(x)\ (n\ge
1).$$ It is well known that (see [B, p.\;151], [G, (3.132)-(3.133)])
$$P_n(x)=\f
1{2^n}\sum_{k=0}^{[n/2]}\b nk(-1)^k\b{2n-2k}nx^{n-2k} =\f 1{2^n\cdot
n!}\cdot\f{d^n}{dx^n}(x^2-1)^n.\tag 1.1$$
 From (1.1)
we see that
$$P_n(-x)=(-1)^nP_n(x), \q P_{2m+1}(0)=0\qtq{and}P_{2m}(0)
=\f{(-1)^m}{2^{2m}}\b{2m}m.\tag 1.2$$
 We also have the following
formula due to Murphy ([G, (3.135)]):
$$P_n(x)=\sum_{k=0}^n\b nk\b{n+k}k\Ls{x-1}2^k=
\sum_{k=0}^n\b{2k}k\b{n+k}{2k}\Ls{x-1}2^k.\tag 1.3$$ We remark that
$\b nk\b{n+k}k=\b{2k}k\b{n+k}{2k}$.
\par Let $\Bbb Z$ be the set of
integers, and for a prime $p$ let $R_p$ be the set of rational
numbers whose denominator is coprime to $p$. Let $[x]$ be the
greatest integer not exceeding $x$, and let $\sls ap$ be the
Legendre symbol. In [S4-S6] the author showed that for any prime
$p>3$ and $t\in R_p$,
$$\align &P_{\f{p-1}2}(t)\e -\Ls{-6}p\sum_{x=0}^{p-1}\Ls{x^3-3(t^2+3)x+2t(t^2-9)}p
\mod p,\tag 1.4
\\&P_{[\f p4]}(t)\e -\Ls{6}p\sum_{x=0}^{p-1}\Ls{x^3-\f{3(3t+5)}2x+9t+7}p\mod
 p,\tag 1.5
\\&P_{[\f  p3]}(t) \e -\Ls
p3\sum_{x=0}^{p-1}\Ls{x^3+3(4t-5)x+2(2t^2-14t+11)}p\mod p.\tag
1.6\endalign$$  In the paper, by using elementary and
straightforward arguments we prove that
$$P_{[\f p6]}(t)
\e -\Ls 3p\sum_{x=0}^{p-1}\Ls{x^3-3x+2t}p\mod p.\tag 1.7 $$
Moreover, for $m,n\in R_p$ with $m\not\e 0\mod p$ we have
$$\sum_{x=0}^{p-1}\Ls{x^3+mx+n}p\e \cases
-(-3m)^{\f{p-1}4}P_{[\f p6]}\ls{3n\sqrt{-3m}}{2m^2}\mod p&\t{if $p\e
1\mod 4$,}
\\-\f{(-3m)^{\f{p+1}4}}{\sqrt{-3m}}P_{[\f p6]}\ls{3n\sqrt{-3m}}{2m^2}\mod
p&\t{if $p\e 3\mod 4$}
\endcases\tag 1.8$$
and
$$\Big(\sum_{x=0}^{p-1}\Big(\frac{x^3+mx+n}p\Big)\Big)^2\equiv
\Big(\frac{-3m}p\Big)
\sum_{k=0}^{[p/6]}\binom{2k}k\binom{3k}k\binom{6k}{3k}
\Big(\frac{4m^3+27n^2}{12^3\cdot 4m^3}\Big)^k\mod p.$$ It is well
known that (see for example [S2, pp.221-222]) the number of points
on the curve $y^2=x^3+mx+n$ over the field $F_p$ with $p$ elements
is given by
$$\#E_p(x^3+mx+n)=p+1+\sum_{x=0}^{p-1}\Ls{x^3+mx+n}p.$$

\par For positive integers $a,b$ and $n$, if
$n=ax^2+by^2$ for some integers $x$ and $y$, we briefly say that
$n=ax^2+by^2$. Recently the author's brother Zhi-Wei Sun[Su1,Su4]
and the author[S4] posed some conjectures for
$\sum_{k=0}^{p-1}\b{2k}k\b{3k}k\b{6k}{3k}/m^k$ modulo ${p^2}$, where
$p>3$ is a prime and $m\in\Bbb Z$ with $p\nmid m$. For example,
Zhi-Wei Sun conjectured that ([Su4, Conjecture 2.8]) for any prime
$p>3$,
$$\sum_{k=0}^{p-1}\f{\b{2k}k\b{3k}k\b{6k}{3k}}{(-96)^{3k}}
\e\cases 0\mod {p^2}&\t{if $\sls p{19}=-1$,}
\\\sls{-6}p(x^2-2p)\mod {p^2}&\t{if $\sls p{19}=1$ and so $4p=x^2+19y^2$.}
\endcases\tag 1.9$$
Using (1.8) and known character sums we determine $P_{[\f
p6]}(x)\mod p$ for $11$ values of $x$ (see Corollaries 2.1-2.11),
and $\sum_{k=0}^{p-1}\b{2k}k\b{3k}k\b{6k}{3k}/m^k\mod {p^2}$ for
$m=-640320^3,$ $-5280^3,-960^3,-96^3,-32^3-15^3,20^3,66^3,255^3,
54000, -12288000$. Thus we solve some conjectures in [Su1,Su4] and
[S4]. For example, we confirm (1.9) in the case $\sls p{19}=-1$ and
prove (1.9) when $\sls p{19}=1$ and the modulus is replaced by $p$.
\par Let $p>3$ be a prime. In the paper we also determine
$\sum_{k=0}^{p-1}\b{3k}k\b{6k}{3k}/864^k\mod{p^2}$ and establish the
 general congruence
$$\sum_{k=0}^{p-1}\b {2k}k\b{3k}k\b{6k}{3k}(x(1-432x))^k\e \Big(
\sum_{k=0}^{p-1}\b {3k}k\b{6k}{3k}x^k\Big)^2\mod {p^2},\tag 1.10$$
and pose some conjectures on supercongruences.
 \subheading{2.
Congruences for $P_{[\f p6]}(x)\mod p$}
 \pro{Lemma 2.1} Let $p$ be an odd prime. Then
 $$\align&\b{\f{p-1}2}k\e \f 1{(-4)^k}\b{2k}k\mod p\qtq{for}
 k=0,1,\ldots,\f{p-1}2,\tag i
  \\&\b{\f{p-1}2-k}{2k}\e \f {\b{6k}{3k}\b{3k}k}{4^{2k}\b{2k}k}\mod p
 \qtq{for}k=0,1,\ldots,\big[\f p3\big],\tag ii
  \\&\b{[\f p3]+k}{2k}\e \f 1{(-27)^k}\b{3k}k\mod p
 \qtq{for $p\not=3$ and}k=0,1,\ldots,\big[\f p3\big].\tag iii\endalign$$
 \endpro
 Proof. For $k\in\{0,1,\ldots,\f{p-1}2\}$ we have $\b{\f{p-1}2}k
 \e \b{-\f 12}k=\f 1{(-4)^k}\b{2k}k\mod p$. Thus (i) holds.
Now suppose $k\in\{0,1,\ldots,[\f p3]\}$. It is clear that
$$\align \b{\f{p-1}2-k}{2k}&=\f{\f{p-1-2k}2\cdot\f{p-3-2k}2\cdots
\f{p-(6k-1)}2}{(2k)!} \e\f {(2k+1)(2k+3)\cdots
(6k-1)}{(-2)^{2k}\cdot (2k)!}
\\&=\f{(6k)!}{4^k(2k)!^2(2(k+1))(2(k+2))\cdots (2(3k))}
=\f{(6k)!}{4^k(2k)!^2\cdot 2^{2k}\cdot \f{(3k)!}{k!}}
\\&=\f {\b{6k}{3k}\b{3k}k}{4^{2k}\b{2k}k}\mod p.
\endalign$$
Thus (ii) is true. (iii) was given by the author in [S4, the proof
of Lemma 2.3].  The proof is now complete.

\pro{Lemma 2.2} Let $p>3$ be a prime and $k\in\{0,1,\ldots,[\f
p{12}]\}$. Then
$$\b{[\f p6]}k\b{2[\f p6]-2k}{[\f p6]}
 \e(-1)^{[\f p6]}3^{3k+\f{1-\sls p3}2}4^{[\f
p6]-k}\b{\f{p-1}2}{[\f p3]-k}\b{\f {p-\sls p3}6+k}{3k+\f{1-\sls
p3}2}\mod p.$$
\endpro
Proof. Using Lemma 2.1(i) we see that
$$\align \b{[\f p6]}k\b{2[\f p6]-2k}{[\f p6]}&=\f{(2[\f
p6]-2k)!}{k!([\f p6]-k)!([\f p6]-2k)!} =\b{2([\f p6]-k)}{[\f
p6]-k}\b{[\f p6]-k}k
\\&\e (-4)^{[\f p6]-k}\b {\f{p-1}2}{[\f p6]-k}\b{[\f p6]-k}k
\\&=(-4)^{[\f p6]-k}\b {\f{p-1}2}k\b{\f{p-1}2-k}{[\f p6]-2k}
\\&\e (-4)^{[\f p6]-2k}\b{2k}k\f{(\f{p-1}2-k)!}{([\f
p6]-2k)!([\f{p+1}3]+k)!}\mod p.
\endalign$$
If $p\e 1\mod 3$, using Lemma 2.1(i) we see that
$$\align \b{\f{p-1}2}{\f{p-1}3-k}\b{\f{p-1}6+k}{3k}&
=\b{\f{p-1}2}{\f{p-1}6+k}\b{\f{p-1}6+k}{3k}=\b{\f{p-1}2}{3k}\b{\f{p-1}2-3k}{\f{p-1}6-2k}
\\&\e\f
1{(-4)^{3k}}\b{6k}{3k}\f{(\f{p-1}2-3k)!}{(\f{p-1}6-2k)!(\f{p-1}3-k)!}
\mod p.\endalign$$ Thus, from the above and Lemma 2.1 we deduce that
$$\align \f{\b{[\f p6]}k\b{2[\f p6]-2k}{[\f
p6]}}{\b{\f{p-1}2}{\f{p-1}3-k}\b{\f{p-1}6+k}{3k}} &\e (-4)^{[\f
p6]-2k+3k}\f{\b{2k}k(\f{p-1}2-k)!(\f{p-1}3-k)!} {\b{6k}{3k}
(\f{p-1}2-3k)!(\f{p-1}3+k)!}
\\&=(-4)^{[\f p6]+k}\f{\b{2k}k\b{\f{p-1}2-k}{2k}}
{\b{6k}{3k}\b{\f{p-1}3+k}{2k}}\e (-4)^{[\f p6]+k}\f{\b{2k}k
\b{6k}{3k}\b{3k}k/(4^{2k}\b{2k}k)} {\b{6k}{3k}\b{3k}k/((-27)^k)}
\\&=(-27)^k(-4)^{[\f p6]-k}=(-1)^{[\f p6]}3^{3k}4^{[\f
p6]-k} \mod p.
\endalign$$
\par If $p\e 2\mod 3$, using Lemma 2.1(i) we see that
$$\align \b{\f{p-1}2}{\f{p-2}3-k}\b{\f{p+1}6+k}{3k+1}&
=\b{\f{p-1}2}{\f{p+1}6+k}\b{\f{p+1}6+k}{3k+1}=\b{\f{p-1}2}{3k+1}\b{\f{p-3}2-3k}
{\f{p-5}6-2k}
\\&\e\f
1{(-4)^{3k+1}}\b{6k+2}{3k+1}\f{(\f{p-3}2-3k)!}{(\f{p-5}6-2k)!(\f{p-2}3-k)!}
\\&\e \f
1{(-4)^{3k}(3k+1)}\b{6k}{3k}\f{(\f{p-1}2-3k)!}{(\f{p-5}6-2k)!(\f{p-2}3-k)!}
\mod p.\endalign$$ Thus, from the above and Lemma 2.1 we deduce that
$$\align \f{\b{[\f p6]}k\b{2[\f p6]-2k}{[\f
p6]}}{\b{\f{p-1}2}{\f{p-2}3-k}\b{\f{p+1}6+k}{3k+1}} &\e 3(-4)^{[\f
p6]-2k+3k}\f{\b{2k}k(\f{p-1}2-k)!(\f{p-2}3-k)!} {\b{6k}{3k}
(\f{p-1}2-3k)!(\f{p-2}3+k)!}
\\&=3(-4)^{[\f p6]+k}\f{\b{2k}k\b{\f{p-1}2-k}{2k}}
{\b{6k}{3k}\b{\f{p-2}3+k}{2k}}\e 3(-4)^{[\f p6]+k}\f{\b{2k}k
\b{6k}{3k}\b{3k}k/(4^{2k}\b{2k}k)} {\b{6k}{3k}\b{3k}k/((-27)^k)}
\\&=3(-27)^k(-4)^{[\f p6]-k} =(-1)^{[\f p6]}3^{3k+1}4^{[\f
p6]-k}\mod p.
\endalign$$
This completes the proof.

\pro{Theorem 2.1} Let $p>3$ be a prime and $m,n\in R_p$ with
$m\not\e 0\mod p$. Then
$$\sum_{x=0}^{p-1}\Ls{x^3+mx+n}p\e \cases
-(-3m)^{\f{p-1}4}P_{[\f p6]}\ls{3n\sqrt{-3m}}{2m^2}\mod p&\t{if $p\e
1\mod 4$,}
\\-\f{(-3m)^{\f{p+1}4}}{\sqrt{-3m}}P_{[\f p6]}\ls{3n\sqrt{-3m}}{2m^2}\mod
p&\t{if $p\e 3\mod 4$.}
\endcases$$
\endpro
Proof. For any positive integer $k$ it is well known that (see [IR,
Lemma 2, p.235])
$$\sum_{x=0}^{p-1}x^k\e \cases p-1\mod p&\t{if $p-1\mid k$,}
\\0\mod p&\t{if $p-1\nmid k$.}
\endcases$$
For $k,r\in\Bbb Z$ with $0\le r\le k\le\f{p-1}2$ we have $0\le
k+2r\le \f{3(p-1)}2$. Thus,
$$\sum_{x=0}^{p-1}x^{k+2r}\e \cases p-1\mod p&\t{if $k=p-1-2r$,}
\\0\mod p&\t{if $k\not=p-1-2r$}
\endcases$$
and therefore
$$\align&\sum_{x=0}^{p-1}(x^3+mx+n)^{\f{p-1}2}\tag 2.1
\\&=\sum_{x=0}^{p-1}\sum_{k=0}^{(p-1)/2}\b{(p-1)/2}k
(x^3+mx)^kn^{\f{p-1}2-k}
\\&=\sum_{x=0}^{p-1}\sum_{k=0}^{(p-1)/2}\b{(p-1)/2}k
\sum_{r=0}^k\b krx^{3r}(mx)^{k-r}n^{\f{p-1}2-k}
\\&=\sum_{r=0}^{(p-1)/2}\sum_{k=r}^{(p-1)/2}
\b{(p-1)/2}k\b krm^{k-r}n^{\f{p-1}2-k}\sum_{x=0}^{p-1}x^{k+2r}
\\&\e (p-1)\sum_{r=0}^{(p-1)/2}\b{(p-1)/2}{p-1-2r}\b{p-1-2r}
rm^{p-1-3r}n^{2r-\f{p-1}2}
\\&=(p-1)\sum_{\f{p-1}4\le r\le \f{p-1}3}\b{(p-1)/2}{p-1-2r}\b{p-1-2r}
rm^{p-1-3r}n^{2r-\f{p-1}2}\mod p.\endalign$$
 If $n\e 0\mod p$, from
the above we deduce that
$$\sum_{x=0}^{p-1}\Ls{x^3+mx+n}p\e
\sum_{x=0}^{p-1}(x^3+mx+n)^{\f{p-1}2}\e \cases
-\b{\f{p-1}2}{\f{p-1}4}m^{\f{p-1}4}\mod p&\t{if $4\mid p-1$,}
\\0\mod p&\t{if $4\mid p-3$.}
\endcases$$
Thus applying (1.2) and Lemma 2.2 (with $k=[\f p{12}]$) we get
$$P_{[\f p6]}(0)=\cases \f 1{(-4)^{[\f p{12}]}}\b{[\f p6]}{[\f
p{12}]}\e (-1)^{[\f p{12}]}3^{\f{p-1}4}\b{\f{p-1}2}{\f{p-1}4}\e
(-3)^{-\f{p-1}4}\b{\f{p-1}2}{\f{p-1}4} \mod p&\t{if $4\mid p-1$,}
\\0&\t{if $4\mid p-3$.}
\endcases$$
 Hence the result is true for $n\e 0\mod p$.
 \par Now we assume $n\not\e 0\mod p$. From (2.1) we see that
$$\align \sum_{x=0}^{p-1}\Ls{x^3+mx+n}p&\e
\sum_{x=0}^{p-1}(x^3+mx+n)^{\f{p-1}2}
\\&\e -\Ls np \sum_{\f{p-1}4\le r\le
\f{p-1}3}\b{(p-1)/2}r\b{\f{p-1}2-r}{p-1-3r}\Ls{n^2}{m^3}^r
\\&= -\Ls
np\sum_{k=0}^{[\f p{12}]}\b{\f{p-1}2}{[\f p3]-k} \b{\f{p-\sls
p3}6+k}{3k+\f{1-\sls p3}2}\Ls{n^2}{m^3}^{[\f p3]-k}\mod
p.\endalign$$
 On the other hand,
$$\aligned &P_{[\f p6]}\sls{3n\sqrt{-3m}}{2m^2}
\\&=2^{-[\f p6]}\sum_{k=0}^{[\f p{12}]}\b{[\f p6]}k(-1)^k\b{2[\f
p6]-2k}{[\f p6]}\Ls{3n\sqrt{-3m}}{2m^2}^{[\f p6]-2k}
\\&=2^{-[\f
p6]}\sum_{k=0}^{[\f p{12}]}\b{[\f p6]}k(-1)^k\b{2[\f p6]-2k}{[\f
p6]}\Ls{3n\sqrt{-3m}}{2m^2}^{\f{p-\sls{-1}p}2}\Big(-\f{27n^2}{4m^3}\Big)^{[\f
p{12}]-k}
\\&=(-1)^{[\f p{12}]}2^{-[\f
p6]}\sum_{k=0}^{[\f p{12}]}\b{[\f p6]}k\b{2[\f p6]-2k}{[\f
p6]}\Ls{3n\sqrt{-3m}}{2m^2}^{\f{p-\sls{-1}p}2}\Big(\f{27n^2}{4m^3}\Big)^{[\f
p3]-k-\f{p-(\f{-1}p)}4}
\\&\e \delta(m,p)^{-1}\Ls np\Ls 3p(-1)^{[\f p{12}]}2^{-[\f
p6]}\sum_{k=0}^{[\f p{12}]}\b{[\f p6]}k\b{2[\f p6]-2k}{[\f
p6]}\Big(\f{27n^2}{4m^3}\Big)^{[\f p3]-k}\mod p,\endaligned$$ where
$$\delta(m,p)=\cases (-3m)^{\f{p-1}4}&\t{if $p\e 1\mod 4$,}
\\\f{(-3m)^{\f{p+1}4}}{\sqrt{-3m}}&\t{if $p\e 3\mod 4$.}
\endcases$$
Hence, by the above and Lemma 2.2 we get
$$\align& \delta(m,p)P_{[\f p6]}\sls{3n\sqrt{-3m}}{2m^2}
\\&\e \Ls np\Ls 3p(-1)^{[\f p{12}]}2^{-[\f p6]}
\\&\q\times\sum_{k=0}^{[\f p{12}]}
(-1)^{[\f p6]}3^{3k+\f{1-\sls p3}2}4^{[\f p6]-k}\b{\f{p-1}2}{[\f
p3]-k}\b{\f {p-\sls p3}6+k}{3k+\f{1-\sls
p3}2}\Big(\f{27n^2}{4m^3}\Big)^{[\f p3]-k}\mod p.\endalign$$ Since
$$\aligned &\Ls 3p(-1)^{[\f p{12}]}2^{-[\f p6]}(-1)^{[\f p6]}3^{3k+\f{1-\sls p3}2}4^{[\f p6]-k}
\Ls {27}4^{[\f p3]-k}
\\&=(-1)^{[\f p{12}]+[\f p6]}\Ls 3p2^{[\f p6]-2[\f p3]}3^{3[\f
p3]+(1-\sls p3)/2}=(-1)^{[\f p{12}]+[\f p6]}\Ls
3p2^{-\f{p-1}2}3^{p-1}
\\&\e (-1)^{[\f p{12}]+[\f p6]}\cdot(-1)^{\f{p-\sls p3}6}\cdot(-1)^{-[\f{p+1}4]}
=(-1)^{2[\f p{12}]}=1\mod p,\endaligned$$ from the above we deduce
that
$$\align \delta(m,p)P_{[\f p6]}\Ls{3n\sqrt{-3m}}{2m^2}
&\e \Ls np\sum_{k=0}^{[\f p{12}]}\b{\f{p-1}2}{[\f p3]-k}\b{\f
{p-\sls p3}6+k}{3k+\f{1-\sls p3}2}\Big(\f{n^2}{m^3}\Big)^{[\f p3]-k}
\\&\e -\sum_{x=0}^{p-1}\Ls{x^3+mx+n}p\mod p.\endalign$$
This completes the proof. \par\q
\newline{\bf Remark 2.1} The congruence (2.1) has been given by the
author in [S6].

\pro{Corollary 2.1} Let $p\not=2,3,11$ be a prime. Then
$$P_{[\f p6]}\Ls{21\sqrt {33}}{121}
\e \cases \ls {33}p (-33)^{\f{p-1}4}2a\mod p&\t{if $4\mid p-1$,
$p=a^2+b^2$ and $4\mid a-1$,} \\0\mod p&\t{if $p\e 3\mod 4$.}
\endcases$$
\endpro
Proof. By [S5, Corollary 2.1 and (2.2)] we have
$$\aligned&\sum_{x=0}^{p-1}\Ls{x^3-11x+14}p
\\&=\Ls 2p\sum_{x=0}^{p-1}\Ls{x^3-4x}p
=\cases (-1)^{\f{p+3}4}2a&\t{if $4\mid p-1$, $p=a^2+b^2$ and $4\mid
a-1$,}
\\0&\t{if $p\e 3\mod 4$.}
\endcases\endaligned\tag 2.2$$
Thus, taking $m=-11$ and $n=14$ in Theorem 2.1 we obtain the result.

 \pro{Corollary 2.2} Let $p>5$ be a prime. Then
$$P_{[\f p6]}\Ls{7\sqrt {10}}{25}\e \cases
(-1)^{\f d2}\sls 5p5^{\f{p-1}4}2c\mod p&\t{if $8\mid p-1$,
$p=c^2+2d^2$ and $4\mid c-1$,}
\\\sls 5p5^{\f {p-3}4}2d\sqrt {10}\mod p&\t{if $8\mid p-3$,
 $p=c^2+2d^2$ and $4\mid d-1$,}\\0\mod p&\t{if $p\e 5,7\mod 8$.}
\endcases$$
\endpro
Proof. Using [S6, Lemma 4.2] we see that
$$\aligned &\sum_{x=0}^{p-1}\Ls{x^3-30x+56}p
\\&=\sum_{x=0}^{p-1}\Ls{(-x)^3-30(-x)+56}p=\Ls{-1}p
\sum_{x=0}^{p-1}\Ls{x^3-30x-56}p
\\&= \cases (-1)^{\f{p+7}8}\sls 3p2c&\t{if $p\e 1\mod 8$, $p=c^2+2d^2$
and $4\mid c-1$,}
\\(-1)^{\f{p-3}8}\sls 3p2c&\t{if $p\e 3\mod 8$,
$p=c^2+2d^2$ and $4\mid c-1$,}
\\0&\t{if $p\e 5,7\mod 8$.}
\endcases\endaligned$$
By [S3, p.1317] we have
$$2^{[\f p4]}\e \cases(-1)^{\f{c^2-1}8}=(-1)^{\f{p-1}8+\f d2}\mod p&\t{if $p=c^2+2d^2\e
1\mod 8$,}
\\(-1)^{\f{c^2-1}8}\f dc=(-1)^{\f{p-3}8}\f dc\mod p&\t{if $p=c^2+2d^2\e 3\mod 8$ with $4\mid c-d$.}
\endcases$$
Now taking $m=-30$ and $n=56$ in Theorem 2.1 and applying the above
we deduce the result.

\pro{Corollary 2.3} Let $p>5$ be a prime. Then
$$P_{[\f p6]}\Ls{11\sqrt 5}{25}\e \cases
 5^{\f{p-1}4}\sls 5p2A\mod p&\t{if $12\mid p-1$, $p=A^2+3B^2$ and
$3\mid A-1$,} \\-5^{\f {p-3}4}\sls 5p2A\sqrt 5\mod p&\t{if $12\mid
p-7$, $p=A^2+3B^2$ and $3\mid A-1$,}\\0\mod p&\t{if $p\e 2\mod 3$.}
\endcases$$
\endpro
Proof. By [S2, Lemma 2.3] (or [S5, Corollary 2.1 (with $t=5/3$) and
(2.3)]) we have
$$\aligned\sum_{x=0}^{p-1}\Ls{x^3-15x+22}p
=\cases -2A&\t{if $3\mid p-1$, $p=A^2+3B^2$ and $3\mid A-1$,}
\\0&\t{if $p\e 2\mod 3$.}
\endcases\endaligned\tag 2.3$$
Thus, taking $m=-15$ and $n=22$ in Theorem 2.1 we obtain the result.
\pro{Corollary 2.4} Let $p>5$ be a prime. Then
$$P_{[\f p6]}\Ls{253\sqrt{10}}{800}
\e\cases -\sls{10}p10^{\f{p-1}4}L\mod p\\\qq\q\qq\t{if $12\mid p-1$,
$4p=L^2+27M^2$ and $3\mid L-1$,}
\\\sls{10}p10^{\f{p-3}4}L\sqrt{10}\mod p\\\qq\qq\q\t{if $12\mid p-7$,
$4p=L^2+27M^2$ and $3\mid L-1$,}
\\0\mod p\ \;\t{if $p\e 2\mod 3$.}
\endcases$$
\endpro
Proof. From [S6, Corollary 3.3] we know that
$$\sum_{x=0}^{p-1}\Ls{x^3-120x+506}p=\cases \sls 2pL&\t{if $3\mid
p-1$, $4p=L^2+27M^2$ and $3\mid L-1$,}\\0&\t{if $p\e 2\mod 3$.}
\endcases\tag 2.4$$
Thus taking $m=-120$ and $n=506$ in Theorem 2.1 we deduce the
result.

\pro{Corollary 2.5} Let $p>7$ be a prime. Then
$$P_{[\f p6]}\Ls{3\sqrt {105}}{25}\e \cases
2\sls p{15} 15^{\f{p-1}4}C\mod p\\\qq\qq\q\t{if $p\e 1,9,25\mod
{28}$, $p=C^2+7D^2$ and $4\mid C-1$,}
\\2\sls p{15}15^{\f {p-3}4}D\sqrt {105}\mod p
\\\qq\qq\q\t{if $p\e 11,15,23\mod {28}$, $p=C^2+7D^2$ and $4\mid D-1$,}\\0\mod p\ \;\t{if
$p\e 3,5,6\mod 7$.}
\endcases$$
\endpro
Proof. Since $(-x-7)^3-35(-x-7)+98=-(x^3+21x^2+112x)$, from [R1,R2]
we see that
$$\aligned &\sum_{x=0}^{p-1}\Ls{x^3-35x+98}p=(-1)^{\f{p-1}2}\sum_{x=0}^{p-1}\Ls{x^3+21
x^2+112x}p
\\&=\cases(-1)^{\f{p+1}2}2\sls C7C&\t{if $p\e 1,2,4\mod 7$ and
 $p=C^2+7D^2$,}\\0&\t{if $p\e 3,5,6\mod
7$.}\endcases\endaligned\tag 2.5$$ Suppose $p\e 1,2,4\mod 7$ and so
$p=C^2+7D^2$. By [S3, p.1317] we have
$$7^{[\f p4]}\e \cases\sls C7\mod p&\t{if $p\e 1,9,25\mod{28}$ and $C\e
1\mod 4$,}
\\-\sls C7\f DC\mod p&\t{if $p\e 11,15,23\mod{28}$ and $D\e
1\mod 4$.}\endcases\tag 2.6$$
 Now taking $m=-35$ and
$n=98$ in Theorem 2.1 and applying all the above we deduce the
result.

\pro{Corollary 2.6} Let $p$ be a prime such that $p\not=2,3,5,7,17$.
\par $(\t{\rm i})$ If $p\e 3,5,6\mod 7$, then $P_{[\f p6]}\sls{171\sqrt{1785}}{85^2}
\e 0\mod p$.
\par $(\t{\rm ii})$ If $p\e 1,2,4\mod 7$ and so $p=C^2+7D^2$ for some $C,D\in\Bbb Z$,
then
$$P_{[\f p6]}\Ls{171\sqrt{1785}}{85^2}
\e\cases \sls{255}p255^{\f{p-1}4}\cdot 2C\mod p&\t{if $4\mid p-1$
and $4\mid C-1$,}
\\\sls{255}p255^{\f{p-3}4}\cdot 2D\sqrt{1785}\mod p&\t{if $4\mid p-3$ and $4\mid D-1$.}
\endcases$$
\endpro
 Proof. From [W, p.296] we know that
$$\aligned&\sum_{x=0}^{p-1}\Ls{(x^2+6x+2)(3x^2+16x)}p\\&=\cases
 -2\sls{-2}p\sls C7C-\sls 3p&\t{if $p\e 1,2,4\mod 7$ and
$p=C^2+7D^2$,}
\\-\sls 3p&\t{if $p\e 3,5,6\mod 7$.}
\endcases\endaligned$$
 As $(x^2+6x+2)(3x^2+16x)=x^4(3+34/x+102/x^2+32/x^3)$, we see
that
$$\align &\sum_{x=0}^{p-1}\Ls{(x^2+6x+2)(3x^2+16x)}p
\\&=\sum_{x=1}^{p-1}\Ls{3+34/x+102/x^2+32/x^3}p=
\sum_{x=1}^{p-1}\Ls{3+34x+102x^2+32x^3}p
\\&=\Ls 2p\sum_{x=1}^{p-1}\Ls{6+68x+204x^2+64x^3}p=\Ls 2p\sum_{x=0}^{p-1}\Ls
{x^3+\f{51}4x^2+17x+6}p-\Ls{12}p\endalign$$ and
$$\align&\sum_{x=0}^{p-1}\Ls
{x^3+\f{51}4x^2+17x+6}p
\\&=\sum_{x=0}^{p-1}\Ls
{(x-\f{17}4)^3+\f{51}4(x-\f{17}4)^2+17(x-\f{17}4)+6}p
=\sum_{x=0}^{p-1}\Ls{x^3-\f{595}{16}x+\f{5586}{64}}p
\\&=\sum_{x=0}^{p-1}\Ls{\sls x4^3-\f{595}{16}\cdot\f
x4+\f{5586}{64}}p=\sum_{x=0}^{p-1}\Ls {x^3-595x+5586}p.
\endalign$$
Now combining all the above we deduce
$$\aligned \sum_{x=0}^{p-1}\Ls {x^3-595x+5586}p=\cases
 (-1)^{\f{p+1}2}2C\sls C7&\t{if $p=C^2+7D^2\e 1,2,4\mod 7$,}
\\0&\t{if $p\e 3,5,6\mod 7$.}
\endcases\endaligned\tag 2.7$$
Taking $m=-595$ and $n=5586$ in Theorem 2.1 and then applying (2.7)
and (2.6) we deduce the result.

\pro{Corollary 2.7} Let $p\not=2,3,11$ be a prime.
\par $(\t{\rm i})$ If $p\e 2,6,7,8,10\mod{11}$, then
$P_{[\f p6]}(\f 7{32}\sqrt{22})\e 0\mod p$.
\par $(\t{\rm ii})$ If $p\e 1,3,4,5,9\mod{11}$ and hence $4p=u^2+11v^2$ for some $u,v\in\Bbb Z$, then
$$P_{[\f p6]}\Big(\f 7{32}\sqrt{22}\Big)\e \cases
-\sls p3(-2)^{\f{p-1}4}u\mod p&\t{if $4\mid p-1$ and $4\mid u-1$,}
\\\sls p32^{\f{p-1}4}u\mod p&\t{if $4\mid p-1$ and $8\mid u-2$,}
\\-\sls p3(-2)^{\f{p-3}4}v\sqrt{22}\mod p&\t{if $4\mid p-3$ and $4\mid v-1$,}
\\\sls p32^{\f{p-3}4}v\sqrt{22}\mod p&\t{if $4\mid p-3$ and $8\mid v-2$.}
\endcases$$
\endpro
Proof. It is known that (see [RP] and [JM])
$$\sum_{x=0}^{p-1}\Ls{x^3-96\cdot 11x+112\cdot 11^2}p=\cases \sls 2p\sls
u{11}u&\t{if $\sls p{11}=1$ and $4p=u^2+11v^2$,}\\0&\t{if $\sls
p{11}=-1$.}\endcases\tag 2.8$$ Thus applying Theorem 2.1 we deduce
$$\align &P_{[\f p6]}(\f 7{32}\sqrt{22})\\&\e\cases -\sls p322^{\f{p-1}4}\sls
u{11}u\mod p&\t{if $\sls p{11}$=1, $4\mid p-1$ and $4p=u^2+11v^2$,}
\\-\sls p322^{\f{p-3}4}\sls
u{11}u\sqrt{22}\mod p&\t{if $\sls p{11}$=1, $4\mid p-3$ and
$4p=u^2+11v^2$,}
\\0\mod p&\t{if $\sls p{11}=-1$.}
\endcases\endalign$$
\par Now assume $\sls p{11}=1$ and so $4p=u^2+11v^2$. If $u\e v\e
1\mod 4$, by [S3, Theorem 4.3] we have
$$(-11)^{[\f p4]}\e\cases \sls u{11}\mod p&\t{if $p\e 1\mod 4$,}
\\\sls u{11}\f vu\mod p&\t{if $p\e 3\mod 4$.}
\endcases$$
If $u\e v\e 0\mod 2$,  by [S3, Corollary 4.6] we have
$$11^{[\f p4]}\e\cases -\sls u{11}\mod p&\t{if $p\e 1\mod 4$ and $8\mid u-2$,}
\\-\sls u{11}\f vu\mod p&\t{if $p\e 3\mod 4$ and $8\mid v-2$.}
\endcases$$
Now combining all the above we derive the result.
\par\q
\par From [RPR], [JM] and [PV] we know that for any prime $p>3$,
$$\aligned &\sum_{x=0}^{p-1}\Ls{x^3-8\cdot 19x+2\cdot 19^2}p
=\cases \sls 2p\sls u{19}u&\t{if $\sls p{19}=1$ and
$4p=u^2+19v^2$,}\\0&\t{if $\sls p{19}=-1$,}
\endcases
\\&\sum_{x=0}^{p-1}\Ls{x^3-80\cdot 43x+42\cdot 43^2}p
=\cases \sls 2p\sls u{43}u&\t{if $\sls p{43}=1$ and
$4p=u^2+43v^2$,}\\0&\t{if $\sls p{43}=-1$,}
\endcases
\\&\sum_{x=0}^{p-1}\Ls{x^3-440\cdot 67x+434\cdot 67^2}p
=\cases \sls 2p\sls u{67}u&\t{if $\sls p{67}=1$ and
$4p=u^2+67v^2$,}\\0&\t{if $\sls p{67}=-1$,}
\endcases\\&\sum_{x=0}^{p-1}\Ls{x^3-80\cdot 23\cdot 29\cdot
163x+14\cdot 11\cdot 19\cdot 127\cdot 163^2}p \\&=\cases \sls 2p\sls
u{163}u&\t{if $\sls p{163}=1$ and $4p=u^2+163v^2$,}\\0&\t{if $\sls
p{163}=-1$.}
\endcases
\endaligned \tag 2.9$$
 Thus, using the method in the proof of
Corollary 2.7 one can similarly determine $P_{[\f p6]}(\f
3{32}\sqrt{114})$, $P_{[\f p6]}(\f {63\sqrt{645}}{1600}),\ P_{[\f
p6]}(\f {651}{96800}\sqrt{22110}),\ P_{[\f p6]}(\f
{557403}{26680^2}\sqrt{1630815})\mod p.$

\pro{Lemma 2.3} Let $p$ be a prime greater than $3$, and let $x$ be
a
 variable. Then
 $$P_{[\f p6]}(x)\e\sum_{k=0}^{[p/6]}\b{6k}{3k}\b{3k}k\Ls{1-x}{864}^k
 \mod p.$$
 \endpro
 Proof. Suppose that $r\in\{1,5\}$ is given by $p\e r\mod 6$. Then clearly
$$\align \b{[\f p6]+k}{2k}&
=\f{(\f{p-r}6+k)(\f{p-r}6+k-1)\cdots(\f{p-r}6-k+1)}{(2k)!}
\\&=\f{(p+6k-r)(p+6k-r-6)\cdots(p-(6k+r-6))}{6^{2k}\cdot (2k)!}
\\&\e (-1)^k\f{(6k-r)(6k-r-6)\cdots(6-r)\cdot r(r+6)\cdots(6k+r-6)}
{6^{2k}\cdot (2k)!}
\\&=\f{(-1)^k\cdot (6k)!}{(2\cdot 4\cdots 6k)(3\cdot 9\cdot 15\cdots (6k-3))\cdot 6^{2k}\cdot (2k)!}
\\&=\f{(-1)^k\cdot (6k)!}{2^{3k}(3k)!\cdot 3^k
\f{(2k)!}{2\cdot 4\cdot 6\cdots 2k}\cdot 36^k(2k)!} \e
\f{(6k)!k!}{(-432)^k(3k)!(2k)!^2}
 \mod p.
\endalign$$
Hence
$$\b{[\f p6]}k\b{[\f p6]+k}k=\b{[\f p6]+k}{2k}\b{2k}k\e
 \f{\b{6k}{3k}\b{3k}k}{(-432)^k}\mod p.\tag 2.10$$
This together with (1.3) yields the result.

\pro{Theorem 2.2} Let $p>3$ be a prime and $m,n\in R_p$ with
$m\not\e 0\mod p$. Then
$$P_{[\f p6]}\Ls n{2m^3}\e \sum_{k=0}^{[p/6]}
\b{6k}{3k}\b{3k}k\Ls{2m^3-n}{12^3m^3}^k\e
-\Ls{3m}p\sum_{x=0}^{p-1}\Ls{x^3-3m^2x+n}p\mod p.$$
\endpro
Proof. Replacing $m$ by $-3m^2$ in Theorem 2.1 and then applying
Lemma 2.3 we deduce the result.\par\q
 \par For positive integers
$a_1,a_2,a_3,a_4$ let
 $$q\prod_{k=1}^{\infty}(1-q^{a_1k})(1-q^{a_2k})(1-q^{a_3k})(1-q^{a_4k})
 =\sum_{n=1}^{\infty}c(a_1,a_2,a_3,a_4;n)q^n\q(|q|<1).$$
For
$(a_1,a_2,a_3,a_4)=(1,1,11,11),(2,2,10,10),(1,3,5,15),(1,2,7,14)$
and $(4,4,8,8)$ it is known that (see [MO, Theorem 1])
$$f(z)=\sum_{n=1}^{\infty}c(a_1,a_2,a_3,a_4;n)q^n \q(q=\t{e}^{2\pi
iz})$$ are weight 2 newforms.
 \pro{Corollary 2.8}  Let $p$ be an odd prime. Then
$$c(1,1,11,11;p)\e P_{[\f p6]}\Ls{19}8\e
(-1)^{\f{p-1}2}\sum_{k=0}^{[p/6]}\f{\b{6k}{3k}\b{3k}k}{256^k}\mod
p.$$\endpro Proof. It is easy to see that the result holds for
$p=3,11$. Now assume $p\not=3,11$. By the well known work of Eichler
in 1954, we have $$|\{(x,y)\in F_p\times F_p:\
y^2+y=x^3-x^2\}|=p-c(1,1,11,11;p).$$
 Since
$$\align &|\{(x,y)\in F_p\times  F_p:\ y^2+y=x^3-x^2\}|
\\&=\Big|\big\{(x,y)\in F_p\times
  F_p:\ \big(y+\f 12\big)^2=x^3-x^2+\f 14\big\}\Big|
\\&=\Big|\big\{(x,y)\in F_p\times  F_p:\
 y^2=x^3-x^2+\f 14\big\}\Big|
\\&=p+\sum_{x=0}^{p-1}\Ls{x^3-x^2+\f 14}p
=p+\sum_{x=0}^{p-1}\Ls{(x+\f 13)^3-(x+\f 13)^2+\f 14}p
\\&=p+\sum_{x=0}^{p-1}\Ls{x^3-\f 13x+\f{19}{108}}p
=p+\sum_{x=0}^{p-1}\Ls {(\f x6)^3-\f 13\cdot \f x6+\f{19}{108}}p
\\&=p+\Ls 6p\sum_{x=0}^{p-1}\Ls{x^3-12x+38}p,
\endalign$$ we obtain
$$c(1,1,11,11;p)=-\Ls 6p\sum_{x=0}^{p-1}\Ls{x^3-12x+38}p.\tag 2.11$$
 Using Theorem 2.2 we see that
$$c(1,1,11,11;p)=-\Ls 6p\sum_{x=0}^{p-1}\Ls{x^3-12x+38}p\e P_{[\f p6]}\Ls{19}8
\mod p.$$ From (1.2) and Lemma 2.3 we have
$$\align P_{[\f p6]}\Ls{19}8&=
(-1)^{[\f p6]}P_{[\f p6]}\Big(-\f {19}8\Big) \e (-1)^{\f{p-1}2}
\sum_{k=0}^{[p/6]}\b{6k}{3k}\b{3k}k\Ls{1+19/8}{864}^k
\\&=(-1)^{\f{p-1}2}\sum_{k=0}^{[p/6]}\f{\b{6k}{3k}\b{3k}k}{256^k}\mod
p.\endalign$$ Thus the result follows.
 \pro{Conjecture 2.1} Let
$p>3$ be a prime. Then
$$\align &c(2,2,10,10;p)=-\Ls p3\sum_{x=0}^{p-1}\Ls{x^3-12x-11}p,
\\&c(2,4,6,12;p)=-\Ls p3\sum_{x=0}^{p-1}\Ls{x^3-39x-70}p,
\\&c(1,3,5,15;p)=-\Ls p3\sum_{x=0}^{p-1}\Ls{x^3-3x-322}p,
\\&c(1,2,7,14;p)=-\Ls p3\sum_{x=0}^{p-1}\Ls{x^3-75x-506}p,
\\&c(4,4,8,8;p)=-\Ls p3\sum_{x=0}^{p-1}\Ls{x^3-99x-378}p.
\endalign$$
\endpro
\par\q If $p>3$ is a prime of the form $4k+3$, from Conjecture 2.1 and [S2, Theorem 2.8] we
deduce that
$$\aligned &p+1-\Ls p3c(2,2,10,10;p)
\\&=p+1+\sum_{x=0}^{p-1}\Ls{x^3-12x-11}p=\#E_p(x^3-12x-11)
\\&=\cases 4N_p-\f{3p-1}2-2\delta(p)&\t{if $p\e 7\mod{12}$,}
\\-4N_p+\f{7p+3}2+2\delta(p)&\t{if $p\e 11\mod{12}$,}
\endcases\endaligned$$
where $N_p$ is the number of $a\in\{0,1,\ldots,p-1\}$ such that
$x^4-4x^2+4x\e a\mod p$ is solvable, and
$$\delta(p)=\cases 0&\t{if $p\e 7,23\mod{40}$,}
\\1&\t{if $p\e 3,27,31,39\mod{40}$,}
\\2&\t{if $p\e 11,19\mod{40}$.}\endcases$$
Hence
$$c(2,2,10,10;p)=\f{5p+1}2-4N_p+2\delta(p)\qtq{for}p\e 3\mod 4.\tag 2.12$$

\pro{Theorem 2.3} Let $p>3$ be a prime, and let $t$ be a variable.
Then
$$P_{[\f p6]}(t)\e -\Ls 3p\sum_{x=0}^{p-1}(x^3-3x+2t)^{\f{p-1}2}\mod
p.\tag 2.13$$
\endpro
Proof. Taking $m=1$ and $n=2t$ in Theorem 2.2 we see that (2.13) is
true for $t=0,1,\ldots,p-1$. Since both sides of (2.13) are
polynomials in $t$ with degree less than $(p-1)/2$, applying
Lagrange's theorem we see that (2.13) holds when $t$ is a variable.

\pro{Theorem 2.4} Let $p>3$ be a prime and let $t$ be a variable.
\par $(\t{\rm i})$ If $t^2+3\not\e 0\mod p$, then
$$P_{\f{p-1}2}(t)\e \cases
(-t^2-3)^{\f{p-1}4}P_{[\f
p6]}\sls{t(t^2-9)\sqrt{t^2+3}}{(t^2+3)^2}\mod p&\t{if $p\e 1\mod
4$,}
\\-\f{(-t^2-3)^{\f{p+1}4}}{\sqrt{t^2+3}}P_{[\f
p6]}\sls{t(t^2-9)\sqrt{t^2+3}}{(t^2+3)^2}\mod p&\t{if $p\e 3\mod
4$.}\endcases$$
\par $(\t{\rm ii})$ If $3t+5\not\e 0\mod p$, then
$$P_{[\f p4]}(t)\e \cases
(6t+10)^{\f{p-1}4}P_{[\f p6]}\sls{(9t+7)\sqrt{6t+10}}{(3t+5)^2}\mod
p&\t{if $p\e 1\mod 4$,}
\\\f{(6t+10)^{\f{p+1}4}}{\sqrt{6t+10}}P_{[\f
p6]}\sls{(9t+7)\sqrt{6t+10}}{(3t+5)^2}\mod p&\t{if $p\e 3\mod
4$.}\endcases$$
\endpro
Proof. Since both sides are polynomials of $t$ with degree at most
$p-3$. It suffices to show that the congruences are true for $t\in
R_p$.  Now combining (1.4)-(1.5) with Theorem 2.1
 we deduce the result.
 \pro{Corollary
2.9} Let $p>3$ be a prime and $m\in R_p$ with $m\not\e 0\mod p$.
Then
$$\align P_{\f {p-1}2}(m^2-3)
&\e \Ls {-2m}pP_{[\f
p6]}\Ls{(m^2-3)(m^2-6)}m\e\sum_{k=0}^{(p-1)/2}\b{2k}k^2\Ls{4-m^2}{32}^k\\&\e
\Ls{2m}p\sum_{k=0}^{[p/6]}\b{6k}{3k}\b{3k}k\Ls{m^4-9m^2-m+18}{864m}^k\mod
p\endalign$$ and
$$\align P_{[\f p4]}\Ls{2m^2-5}3
&\e  \Ls {2m}pP_{[\f
p6]}\Ls{3m^2-4}{m^3}\e\sum_{k=0}^{[p/4]}\b{4k}{2k}\b{2k}k\Ls{4-m^2}{192}^k\\&\e
\Ls{2m}p\sum_{k=0}^{[p/6]}\b{6k}{3k}\b{3k}k\Ls{(m+1)(m-2)^2}{864m^3}^k\mod
p.\endalign$$
\endpro
Proof. Taking $t=m^2-3$ in Theorem 2.4(i) and then applying [S4,
(2.4)] and Lemma 2.3 we deduce the first congruence.
 Taking $t=(2m^2-5)/3$ in Theorem 2.4(ii) and then applying
[S5, Theorem 2.1(ii)] and Lemma 2.3 we deduce the second congruence.

 \pro{Theorem 2.5} Let $p>3$ be a prime and let $t$
be a variable. Then
$$P_{[\f p3]}(t)\e \cases (5-4t)^{\f{p-1}4}P_{[\f p6]}\big(
\f{2t^2-14t+11}{(5-4t)^2}\sqrt{5-4t}\big)\mod p&\t{if $p\e 1\mod
4$,}
\\-\f{(5-4t)^{\f{p+1}4}}{\sqrt{5-4t}}P_{[\f p6]}\big(
\f{2t^2-14t+11}{(5-4t)^2}\sqrt{5-4t}\big)\mod p&\t{if $p\e 3\mod
4$.}
\endcases$$
\endpro
Proof. Since both sides are polynomials in $t$ with degree at most
$p-2$. It suffices to show that the congruence is true for all $t\in
R_p$ with $t\not\e \f 54\mod p$. Set $m=3(4t-5)$ and
$n=2(2t^2-14t+11)$. Then
$$\f{3n\sqrt{-3m}}{2m^2}=\f{(2t^2-14t+11)\sqrt{5-4t}}{(5-4t)^2}.$$
Thus, by (1.6) and Theorem 2.1 we have
$$\aligned P_{[\f p3]}(t)&\e-\Ls p3\sum_{x=0}^{p-1}\Ls{x^3+mx+n}p
\\&\e\cases \ls p3(9(5t-4))^{\f{p-1}4}P_{[\f p6]}\big(
\f{2t^2-14t+11}{(5-4t)^2}\sqrt{5-4t}\big)\mod p&\t{if $4\mid p-1$,}
\\\ls p3\f{(9(5t-4))^{\f{p+1}4}}{\sqrt {9(5-4t)}}P_{[\f p6]}\big(
\f{2t^2-14t+11}{(5-4t)^2}\sqrt{5-4t}\big)\mod p&\t{if $4\mid
p-3$.}\endcases\endaligned$$ For $p\e 1\mod 4$ we have
$9^{\f{p-1}4}\ls p3 \e \ls 3p\ls p3=1\mod p$, For $p\e 3\mod 4$ we
have $9^{\f{p+1}4}\f 13\ls p3 \e \ls 3p\ls p3=-1\mod p$. Thus the
result follows.
 \pro{Corollary 2.10} Let $p>3$ be a prime and
$m\in R_p$ with $m\not\e 0\mod p$. Then
$$\align P_{[\f p3]}\Ls{5-m^2}4&
\e\Ls mpP_{[\f p6]}\Ls{m^4+18m^2-27}{8m^3}\e
\sum_{k=0}^{[p/3]}\b{2k}k\b{3k}k\Ls{m^2-1}{216}^k
 \\& \e
\Ls{-m}p\sum_{k=0}^{[p/6]}\b{6k}{3k}\b{3k}k\Ls{(m+1)(3-m)^3}
{2^8\cdot 3^3m^3}^k\mod p.\endalign$$
\endpro
Proof. Taking $t=\f{5-m^2}4$ in Theorem 2.5 and then applying [S4,
Lemma 2.3] and Lemma 2.3 we deduce the result.

\pro{Corollary 2.11} Let $p>3$ be a prime. Then
$$\aligned&P_{[\f p3]}\Ls{7\pm 3\sqrt 3}2\\&\e\cases 2a(-3\pm 2\sqrt
3)^{\f{p-1}4}\mod p&\t{if $p\e 1\mod 4$ and $p=a^2+b^2$ with $a\e
1\mod 4$,}
\\0\mod p&\t{if $p\e 3\mod 4$.}
\endcases\endaligned$$
\endpro
Proof. Set $t=(7\pm 3\sqrt 3)/2$. Then $2t^2-14t+11=0$. Thus, from
Theorem 2.5 and the congruence for $P_{[\f p6]}(0)$ in the proof of
Theorem 2.1 we deduce
$$P_{[\f p3]}\Ls{7\pm 3\sqrt 3}2\e\cases (-9\mp 6\sqrt
3)^{\f{p-1}4}\e (-3\pm 2\sqrt
3)^{\f{p-1}4}\b{\f{p-1}2}{\f{p-1}4}\mod p&\t{if $p\e 1\mod 4$,}
\\0\mod p&\t{if $p\e 3\mod 4$.}
\endcases$$
It is well known that $\b{\f{p-1}2}{\f{p-1}4}\e 2a\mod p$ for $p\e
1\mod 4$ (see [BEW, p.269]). Thus the corollary is proved.

\pro{Theorem 2.6} Let $p>3$ be a prime and $m,n\in R_p$ with
$m\not\e 0\mod p$. Then
$$\aligned &\sum_{x=0}^{p-1}\Ls{x^3+mx+n}p
\\&\e\cases  -(-3m)^{\f{p-1}4}
\sum_{k=0}^{[p/12]}\b{[\f p{12}]}k\b{[\f {5p}{12}]}k
(\f{4m^3+27n^2}{4m^3})^k\mod p&\t{if $4\mid p-1$,}
\\-\f{3n}{2m^2}(-3m)^{\f{p+1}4}
\sum_{k=0}^{[p/12]}\b{[\f p{12}]}k\b{[\f {5p}{12}]}k
(\f{4m^3+27n^2}{4m^3})^k\mod p&\t{if $4\mid p-3$.}
\endcases
\endaligned$$
\endpro
Proof. Let $P_n^{(\alpha,\beta)}(x)$ be the Jacobi polynomial
defined by
$$P_n^{(\alpha,\beta)}(x)=\f 1{2^n}\sum_{k=0}^n\b{n+\alpha}k\b{n+\beta}{n-k}
(x+1)^k(x-1)^{n-k}.$$ It is known that (see [AAR, p.315])
$$P_{2n}(x)=P_n^{(0,-\f 12)}(2x^2-1)\qtq{and}
P_{2n+1}(x)=xP_n^{(0,\f 12)}(2x^2-1).\tag 2.14$$ From [B, p.170] we
know that
$$\align P_n^{(\alpha,\beta)}(x)&=\b{n+\alpha}n\sum_{k=0}^{\infty}
\f{(-n)_k(n+\alpha+\beta+1)_k}{(\alpha+1)_k\cdot k!}\Ls{1-x}2^k
\\&=\b{n+\alpha}n\sum_{k=0}^n\f{\b
nk\b{-n-\alpha-\beta-1}k}{\b{-1-\alpha}k}\Ls{x-1}2^k.\endalign$$
Thus,
$$P_n^{(0,\beta)}(x)=\sum_{k=0}^n\b nk\b{-n-\beta-1}k\Ls{1-x}2^k.
\tag 2.15$$ Hence, if $p\e 1\mod 4$, then $[\f p6]=2[\f p{12}]$ and
so
$$\align P_{[\f
p6]}\Ls{3n\sqrt{-3m}}{2m^2}&=P_{[\f p{12}]}^{(0,-\f 12)}\Big(2\cdot
\f{-27n^2}{4m^3}-1\Big) =\sum_{k=0}^{[\f p{12}]}\b{[\f p{12}]}k
\b{-\f 12-[\f p{12}]}k\Big(1-\f{-27n^2}{4m^3}\Big)^k
\\&\e \sum_{k=0}^{[\f p{12}]}\b{[\f p{12}]}k
\b{\f {p-1}2-[\f p{12}]}k\Ls{4m^3+27n^2}{4m^3}^k
\\&=\sum_{k=0}^{[\f p{12}]}\b{[\f p{12}]}k
\b{[\f {5p}{12}]}k\Ls{4m^3+27n^2}{4m^3}^k\mod p;\endalign$$ if $p\e
3\mod 4$, then $[\f p6]=2[\f p{12}]+1$ and so
$$\align P_{[\f
p6]}\Ls{3n\sqrt{-3m}}{2m^2}&=\f{3n\sqrt{-3m}}{2m^2}P_{2[\f
p{12}]+1}^{(0,\f 12)}\Big(2\cdot \f{-27n^2}{4m^3}-1\Big)
\\&=\f{3n\sqrt{-3m}}{2m^2}\sum_{k=0}^{[\f p{12}]}\b{[\f p{12}]}k \b{-\f
32-[\f p{12}]}k\Big(1-\f{-27n^2}{4m^3}\Big)^k
\\&\e \f{3n\sqrt{-3m}}{2m^2}\sum_{k=0}^{[\f p{12}]}\b{[\f p{12}]}k
\b{\f {p-3}2-[\f p{12}]}k\Ls{4m^3+27n^2}{4m^3}^k
\\&=\f{3n\sqrt{-3m}}{2m^2}\sum_{k=0}^{[\f p{12}]}\b{[\f p{12}]}k
\b{[\f {5p}{12}]}k\Ls{4m^3+27n^2}{4m^3}^k\mod p.\endalign$$ Now
combining the above with Theorem 2.1 we deduce the result.

 \subheading{3. A
general congruence modulo $p^2$}
 \pro{Lemma 3.1} For any nonnegative integer $n$ we have
$$\sum_{k=0}^n\b {2k}k\b{3k}k\b{6k}{3k}\b k{n-k}(-432)^{n-k}
=\sum_{k=0}^n\b {3k}k\b{6k}{3k}\b{3(n-k)}{n-k}\b{6(n-k)}{3(n-k)}.$$
\endpro
Proof. Let $m$ be a nonnegative integer. For $k\in\{0,1,\ldots,m\}$
set
$$\align &F_1(m,k)=\b {2k}k\b{3k}k\b{6k}{3k}\b k{m-k}(-432)^{m-k},
\\&F_2(m,k)=\b {3k}k\b{6k}{3k}\b{3(m-k)}{m-k}\b{6(m-k)}{3(m-k)}.
\endalign$$
For $k\in\{0,1,\ldots,m+1\}$ set
$$\align &G_1(m,k)=\f{186624k^2(m+2)(m+1-2k)}{k+1}
\b{2k}k\b{3k}k\b{6k}{3k}\b{k+1}{m+2-k}(-432)^{m-k},
\\&G_2(m,k)=
\frac{12k^2(36m^2-36mk+129m-62k+114)}{(m-k+2)^2}
\\&\qq\qq\q\times\b{3k}k\b{6k}{3k}\b{3(m-k+1)}{m-k+1}\b{6(m-k+1)}
{3(m-k+1)}.
\endalign$$ For
$i=1,2$ and $k\in\{0,1,\ldots,m\}$, it is easy to check that
$$\aligned &(m+2)^3F_i(m+2,k)-24(2m+3)(18m^2+54m+41)F_i(m+1,k)
\\&\qq + 20736(m+1)(3m+1)(3m+5)F_i(m,k)=G_i(m,k+1)-G_i(m,k).
\endaligned\tag 3.1$$
Set $S_i(n)=\sum_{k=0}^nF_i(n,k)$ for $n=0,1,2,\ldots$. Then
$$\align &(m+2)^3(S_i(m+2)-F_i(m+2,m+2)-F_i(m+2,m+1))
\\&\q-24(2m+3)(18m^2+54m+41) (S_i(m+1)-F_i(m+1,m+1))
\\&\q+ 20736(m+1)(3m+1)(3m+5)S_i(m)
\\&=(m+2)^3\sum_{k=0}^mF_i(m+2,k)-24(2m+3)(18m^2+54m+41)
\sum_{k=0}^mF_i(m+1,k)
\\&\q+ 20736(m+1)(3m+1)(3m+5)\sum_{k=0}^mF_i(m,k)
\\&=\sum_{k=0}^m (G_i(m,k+1)-G_i(m,k))
=G_i(m,m+1)-G_i(m,0) =G_i(m,m+1).\endalign$$ Thus, for $i=1,2$ and
$m=0,1,2,\ldots$,
$$\aligned &(m+2)^3S_i(m+2)-24(2m+3)(18m^2+54m+41)S_i(m+1)
\\&\q+ 20736(m+1)(3m+1)(3m+5)S_i(m)
\\&=G_i(m,m+1)+(m+2)^3(F_i(m+2,m+2)+F_i(m+2,m+1))
\\&\q-24(2m+3)(18m^2+54m+41) F_i(m+1,m+1)=0.\endaligned\tag 3.2$$
Since $S_1(0)=1=S_2(0)$ and $S_1(1)=120=S_2(1)$, from (3.2) we
deduce $S_1(n)=S_2(n)$ for all $n=0,1,2,\ldots$. This completes the
proof.

\par For given prime $p$ and integer $n$, if $p^{\alpha}\mid n$ but
$p^{\alpha+1}\nmid n$, we say that $p^{\alpha}\ \Vert\ n$.

\pro{Lemma 3.2} Let $p$ be an odd prime and
$k,r\in\{0,1,\ldots,p-1\}$ with $k+r\ge p$. Then
$$\b{3k}k\b{6k}{3k}\b{3r}r\b{6r}{3r}\e 0\mod {p^2}.$$
\endpro
Proof. If $k>\f{5p}6$, then $p^5\mid (6k)!$, $p\ \Vert\ (2k)!$,
$p^2\ \Vert\ (3k)!$ and so
$\b{3k}k\b{6k}{3k}=\f{(6k)!}{k!(2k)!(3k)!}\e 0\mod {p^2}$. If
$\f{2p}3\le k<\f{5p}6$, then $2p\le 3k<3p$, $4p\le 6k<5p$, $p^4\
\Vert\ (6k)!$, $p^2\ \Vert\ (3k)!$, $p\ \Vert\ (2k)!$ and so
$\b{3k}k\b{6k}{3k}=\f{(6k)!}{k!(2k)!(3k)!}\e 0\mod p$. If $\f
p2<k<\f{2p}3$, then $p<3k<2p$, $3p<6k<4p$, $p^3\mid (6k)!$, $p\
\Vert\ (2k)!$, $p\ \Vert\ (3k)!$ and so
$\b{3k}k\b{6k}{3k}=\f{(6k)!}{k!(2k)!(3k)!}\e 0\mod p$. If $\f p3\le
k<\f p2$, then $2k<p\le 3k<2p\le 6k$, $p^2\mid (6k)!$, $p\nmid
(2k)!$, $p\ \Vert\ (3k)!$ and so
$\b{3k}k\b{6k}{3k}=\f{(6k)!}{k!(2k)!(3k)!}\e 0\mod p$. If $\f
p6<k<\f p3$, then $3k<p<6k$ and so
$\b{3k}k\b{6k}{3k}=\f{(6k)!}{k!(2k)!(3k)!}\e 0\mod p$.
\par From the above we see that $p\mid \b{3k}k\b{6k}{3k}$ for $k>\f
p6$. Therefore, if $k>\f p6$ and $r>\f p6$, then $\b{3k}k\b{6k}{3k}
\b{3r}r\b{6r}{3r}\e 0\mod{p^2}$. If $r<\f p6$, then $k\ge
p-r>\f{5p}6$ and so $p^2\mid \b{3k}k\b{6k}{3k}$ by the above. If
$k<\f p6$, then $r\ge p-k>\f{5p}6$ and so $p^2\mid
\b{3r}r\b{6r}{3r}$ by the above.
\par Now putting all the above together we prove the lemma.

 \pro{Theorem 3.1} Let $p$ be an odd prime and let  $x$ be a
variable. Then
$$\sum_{k=0}^{p-1}\b {2k}k\b{3k}k\b{6k}{3k}(x(1-432x))^k\e \Big(
\sum_{k=0}^{p-1}\b {3k}k\b{6k}{3k}x^k\Big)^2\mod {p^2}.$$
\endpro
Proof. It is clear that
$$\align &\sum_{k=0}^{p-1}\b {2k}k\b{3k}k\b{6k}{3k}(x(1-432x))^k
\\&=\sum_{k=0}^{p-1}\b {2k}k\b{3k}k\b{6k}{3k}x^k\sum_{r=0}^k\b kr(-432x)^r
\\&=\sum_{m=0}^{2(p-1)}x^m\sum_{k=0}^{min\{m,p-1\}}\b {2k}k\b{3k}k\b{6k}{3k}
\b k{m-k}(-432)^{m-k}.\endalign$$
 Suppose $p\le m\le 2p-2$ and $0\le
k\le p-1$. If $k\ge \f {2p}3$, then $2p\le 3k<3p$, $6k\ge 4p$,
$p^3\nmid (3k)!$, $p^4\mid (6k)!$ and so $\b
{2k}k\b{3k}k\b{6k}{3k}=\f{(6k)!}{(3k)!k!^3}\e 0\mod {p^2}$. If $\f
p2<k<\f{2p}3$, then $3k<2p$, $6k>3p$, $p^2\nmid (3k)!$ and $p^3\mid
(6k)!$ and so $\b {2k}k\b{3k}k\b{6k}{3k}=\f{(6k)!}{(3k)!k!^3}\e
0\mod {p^2}$.  If $k<\f p2$, then $m-k\ge p-k>k$ and so $\b
k{m-k}=0$. Thus, from the above and Lemma 3.1 we deduce that
$$\align &\sum_{k=0}^{p-1}\b {2k}k\b{3k}k\b{6k}{3k}(x(1-432x))^k
\\&\e \sum_{m=0}^{p-1}x^m\sum_{k=0}^m\b {2k}k\b{3k}k\b{6k}{3k}
\b k{m-k}(-432)^{m-k}
\\&=\sum_{m=0}^{p-1}x^m\sum_{k=0}^m\b {3k}k\b{6k}{3k}
\b{3(m-k)}{m-k}\b{6(m-k)}{3(m-k)}
\\&=\sum_{k=0}^{p-1}\b{3k}k\b{6k}{3k}x^k\sum_{m=k}^{p-1}
\b{3(m-k)}{m-k}\b{6(m-k)}{3(m-k)}x^{m-k}
\\&=\sum_{k=0}^{p-1}\b{3k}k\b{6k}{3k}x^k
\sum_{r=0}^{p-1-k}\b{3r}r\b{6r}{3r}x^r
\\&=\sum_{k=0}^{p-1}\b{3k}k\b{6k}{3k}x^k\Big(\sum_{r=0}^{p-1}\b{3r}r
\b{6r}{3r}x^r -\sum_{r=p-k}^{p-1}\b{3r}r\b{6r}{3r}x^r\Big)
\\&=\Big(\sum_{k=0}^{p-1}\b{3k}k\b{6k}{3k}x^k\Big)^2
-\sum_{k=0}^{p-1}\b{3k}k\b{6k}{3k}x^k
\sum_{r=p-k}^{p-1}\b{3r}r\b{6r}{3r}x^r
 \mod{p^2}.\endalign$$
By Lemma 3.2 we have $p^2\mid \b{3k}k\b{6k}{3k}\b{3r}r\b{6r}{3r}$
for $0\le k\le p-1$ and $p-k\le r\le p-1$. Thus
$$\sum_{k=0}^{p-1}\b{3k}k\b{6k}{3k}x^k
\sum_{r=p-k}^{p-1}\b{3r}r\b{6r}{3r}x^r
 \e 0\mod{p^2}.$$
Now combining all the above we obtain the result. \pro{Corollary
3.1} Let $p$ be a prime greater than $3$ and $m\in R_p$ with
$m\not\e 0\mod p$. Then
$$\sum_{k=0}^{p-1}\f{\b{2k}k\b{3k}k\b{6k}{3k}}{m^k}
\e\Big(\sum_{k=0}^{p-1}\b{3k}k\b{6k}{3k}
\Ls{1-\sqrt{1-1728/m}}{864}^k\Big)^2\mod
{p^2}.$$
\endpro
Proof. Taking $x=\f{1-\sqrt{1-1728/m}}{864}$ in Theorem 3.1 we
deduce the result.
 \pro{Lemma 3.3} Let $p$
be a prime of the form $4k+1$ and $p=a^2+b^2(a,b\in\Bbb Z)$ with
$a\e 1\mod 4$. Then
$$P_{[\f p6]}(0)\e \b{\f{p-1}2}{[\f p{12}]}\e
 \cases 2a&\t{if $p\e 1\mod{12}$ and $3\nmid a$,}\\-2a&\t{if $p\e
1\mod{12}$ and $3\mid a$,}
\\2b&\t{if $p\e 5\mod{12}$ and $3\mid a-b$.}
\endcases$$
\endpro
Proof. By Lemma 2.1(i) and  the proof of Theorem 2.1  we have
$$P_{[\f p6]}(0)=\f 1{(-4)^{[\f p{12}]}}\b{[\f p6]}{[\f p{12}]}
\e \b{\f{p-1}2}{[\f p{12}]}\e
 (-3)^{-\f{p-1}4}\b{\f{p-1}2}{\f{p-1}4}\mod p.$$ By Gauss'
congruence ([BEW, p.269]), $\b{\f{p-1}2}{\f{p-1}4}\e 2a\mod p$ . By
[S1, Theorem 2.2],
$$(-3)^{\f{p-1}4}\e \cases 1\mod p&\t{if $p\e 1\mod {12}$ and $3\nmid
a$,}\\-1\mod p&\t{if $p\e 1\mod {12}$ and $3\mid a$,}
  \\-\f ba\e\f ab\mod p&\t{if $p\e 5\mod{12}$ and $b\e a\mod 3$.}
\endcases$$
   Thus the result follows.

\par Let $p>3$ be a prime. By the work of Mortenson[M] and
Zhi-Wei Sun[Su3],
$$\sum_{k=0}^{p-1}\f{\b{2k}k\b{3k}k\b{6k}{3k}}{1728^k}\e \cases\sls
p3(4a^2-2p)\mod{p^2}&\t{if $p=a^2+b^2\e 1\mod 4$ and $2\nmid a$,}
\\0\mod{p^2}&\t{if $p\e 3\mod 4$.}\endcases\tag 3.3$$ In [Su1]
Zhi-Wei Sun conjectured that
$$\sum_{k=0}^{p-1}\f{\b{6k}{3k}\b{3k}k}{864^k}\e \cases 0\mod
{p^2}&\t{if $p\e 7,11\mod {12}$,}
\\(-1)^{[\f a6]}(2a-\f p{2a})\mod {p^2}&\t{if $12\mid p-1$, $p=a^2+b^2$ and $4\mid a-1$,}
\\\ls{ab}3(2b-\f p{2b}) \mod {p^2}&\t{if $12\mid p-5$, $p=a^2+b^2$ and $4\mid a-1$.}
\endcases$$ In [Su2], Zhi-Wei Sun confirmed the conjecture in the
case $p\e 3\mod 4$.
\par Now we prove the above conjecture for primes $p\e 1\mod 4$.

\pro{Theorem 3.2} Let $p$ be a prime of the form $4k+1$ and so
$p=a^2+b^2$ with $a,b\in\Bbb Z$ and $4\mid a-1$. Then
$$\sum_{k=0}^{p-1}\f{\b{6k}{3k}\b{3k}k}{864^k}
\e \cases 2a-\f p{2a}\mod{p^2}&\t{if $p\e 1\mod{12}$ and $3\nmid
a$,}\\-2a+\f p{2a}\mod{p^2}&\t{if $p\e 1\mod{12}$ and $3\mid a$,}
\\2b-\f p{2b}\mod{p^2}&\t{if $p\e 5\mod{12}$ and $3\mid a-b$.}
\endcases$$
\endpro
Proof. From Lemma 3.3 we have $ P_{[\f p6]}(0) \e 2r\mod p,$
 where
$$r= \cases a&\t{if $p\e 1\mod{12}$ and $3\nmid
a$,}\\-a&\t{if $p\e 1\mod{12}$ and $3\mid a$,}
\\b&\t{if $p\e 5\mod{12}$ and $3\mid a-b$.}
\endcases$$
By the proof of Lemma 3.2 we have $p\mid\b{6k}{3k}\b{3k}k$ for
$p>k>\f p6$. Thus, applying Lemma 2.3 and the above we get
$$\sum_{k=0}^{p-1}\f{\b{6k}{3k}\b{3k}k}{864^k}
\e \sum_{k=0}^{[p/6]}\f{\b{6k}{3k}\b{3k}k}{864^k} \e P_{[\f
p6]}(0)\e 2r\mod p.$$ Set
$\sum_{k=0}^{p-1}\f{\b{6k}{3k}\b{3k}k}{864^k}=2r+qp$. Using
Corollary 3.1 we see that
$$\sum_{k=0}^{p-1}\f{\b{2k}k\b{3k}k\b{6k}{3k}}{1728^k}
\e\Big(\sum_{k=0}^{p-1}\f{\b{6k}{3k}\b{3k}k}{864^k}\Big)^2
=(2r+qp)^2\e 4r^2+4rqp\mod{p^2}.$$  Thus, applying (3.3) we obtain
$\sls p3(4a^2-2p)\e 4r^2+4rqp\mod{p^2}.$ Hence $q\e -\f 1{2r}\mod p$
and the proof is complete.
 \subheading{4. Congruences for
$\sum_{k=0}^{p-1}\b{2k}k\b{3k}k\b{6k}{3k}/m^k$}

\pro{Theorem 4.1} Let $p>3$ be a prime, $m\in R_p$, $m\not\e 0\mod
p$ and $t=\sqrt{1-1728/m}$. Then
$$\sum_{k=0}^{p-1}\f{\b{2k}k\b{3k}k\b{6k}{3k}}{m^k}\e
P_{[\f p6]}(t)^2\e
\Big(\sum_{x=0}^{p-1}(x^3-3x+2t)^{\f{p-1}2}\Big)^2\mod p.$$
Moreover, if $P_{[\f p6]}(t)\e 0\mod p$ or
$\sum_{x=0}^{p-1}(x^3-3x+2t)^{\f{p-1}2}\e 0\mod p$, then
$$\sum_{k=0}^{p-1}\f{\b{2k}k\b{3k}k\b{6k}{3k}}{m^k}\e 0\mod{p^2}.$$
\endpro
Proof. Since $\f{1-t}{864}(1-432\cdot \f{1-t}{864})=\f 1m$, by
Theorem 3.1 we have
$$\sum_{k=0}^{p-1}\f{\b{2k}k\b{3k}k\b{6k}{3k}}{m^k}\e
\Big(\sum_{k=0}^{p-1}\b{3k}k\b{6k}{3k}\Ls{1-t}{864}^k\Big)^2\mod{p^2}.
\tag 4.1$$ From the proof of Lemma 3.2 we know that $p\mid
\b{3k}k\b{6k}{3k}$ for $[\f p6]<k<p$. Thus, using Lemma 2.3 and
Theorem 2.3 we see that
$$\sum_{k=0}^{p-1}\b{2k}k\b{3k}k\Ls{1-t}{864}^k
\e P_{[\f p6]}(t)\e -\Ls
p3\sum_{x=0}^{p-1}(x^3-3x+2t)^{\f{p-1}2}\mod p.$$ This together with
(4.1) yields the result.

 \pro{Theorem 4.2} Let $p>3$ be a prime
and $m,n\in R_p$ with $m\not\e 0\mod p$. Then
$$\align \Big(\sum_{x=0}^{p-1}\Ls{x^3+mx+n}p\Big)^2&
\e
\Ls{-3m}p\sum_{k=0}^{p-1}\b{2k}k\b{3k}k\b{6k}{3k}\Ls{4m^3+27n^2}{12^3\cdot
4m^3}^k\\&\e
\Ls{-3m}p\sum_{k=0}^{[p/6]}\b{2k}k\b{3k}k\b{6k}{3k}\Ls{4m^3+27n^2}{12^3\cdot
4m^3}^k\mod p.\endalign$$ Moreover, if
$\sum_{x=0}^{p-1}\sls{x^3+mx+n}p=0$, then
$$\sum_{k=0}^{p-1}\b{2k}k\b{3k}k\b{6k}{3k}\Ls{4m^3+27n^2}{12^3\cdot
4m^3}^k\e 0\mod{p^2}.$$
\endpro
Proof. By the proof of Lemma 3.2 we have $p\mid \b{3k}k\b{6k}{3k}$
for $\f p6<k<p$. We first assume $4m^3+27n^2\e 0\mod p$. As
$x^3+mx+n\e (x-\f{3n}m)(x+\f{3n}{2m})^2\mod p$ we see that
$$\align\sum_{x=0}^{p-1}\Ls{x^3+mx+n}p&=\sum_{x=0}^{p-1}\Ls{(x-\f{3n}m)(x+\f{3n}{2m})^2}p
=\sum\Sb x=0\\x\not\e -\f{3n}{2m}\mod p\endSb^{p-1}\Ls{x-\f{3n}m}p
\\&=\sum_{t=0}^{p-1}\Ls
tp-\Ls{-\f{3n}{2m}-\f{3n}m}p=-\Ls{-2mn}p.\endalign$$
 Since $m\not\e 0\mod
p$ we have $n\not\e 0\mod p$ and so
$\sum_{x=0}^{p-1}\sls{x^3+mx+n}p=-\sls{-2mn}p\not=0$.
 Thus
the result holds in this case.
\par Now we assume $4m^3+27n^2\not\e
0\mod p$. Set $t=\f{3n\sqrt{-3m}}{2m^2}$ and $m_1=\f{1728\cdot
4m^3}{4m^3+27n^2}$. Then $t=\sqrt{1-\f{1728}{m_1}}$. From Theorems
2.1 and 4.1 we have
$$ \Big(\sum_{x=0}^{p-1}\Ls{x^3+mx+n}p\Big)^2\e
(-3m)^{\f{p-1}2}P_{[\f p6]}(t)^2\e \Ls{-3m}p\sum_{k=0}^{p-1}
\f{\b{2k}k\b{3k}k\b{6k}{3k}}{m_1^k} \mod p.$$ If
$\sum_{x=0}^{p-1}\sls{x^3+mx+n}p=0$, using Theorems 2.1 and 4.1 we
see that $P_{[\f p6]}(t)\e 0\mod p$ and so $\sum_{k=0}^{p-1}
\f{\b{2k}k\b{3k}k\b{6k}{3k}}{m_1^k}\e 0\mod{p^2}$.
 This completes the proof.

\pro{Theorem 4.3 ([Su4, Conjecture 2.7])} Let $p\not=2,11$ be a
prime. Then
$$\sum_{k=0}^{p-1}\f{\b{2k}k\b{3k}k\b{6k}{3k}}{(-32)^{3k}}
\e\cases 0\mod {p^2}&\t{if $\sls p{11}=-1$,}
\\\sls{-2}px^2\mod p&\t{if $\sls p{11}=1$ and so $4p=x^2+11y^2$.}
\endcases$$
\endpro
Proof. Taking $m=-96\cdot 11$ and $n=112\cdot 11^2$ in Theorem 4.2
and then applying (2.8) we deduce the result.

\pro{Theorem 4.4 ([Su4, Conjecture 2.8])} Let $p\not=2,3,19$ be a
prime. Then
$$\sum_{k=0}^{p-1}\f{\b{2k}k\b{3k}k\b{6k}{3k}}{(-96)^{3k}}
\e\cases 0\mod {p^2}&\t{if $\sls p{19}=-1$,}
\\\sls{-6}px^2\mod p&\t{if $\sls p{19}=1$ and so $4p=x^2+19y^2$.}
\endcases$$
\endpro
Proof. Taking $m=-8\cdot 19$ and $n=2\cdot 19^2$ in Theorem 4.2 and
then applying (2.9) we deduce the result.

\pro{Theorem 4.5 ([Su4, Conjecture 2.9])} Let $p\not=2,3,5,43$ be a
prime. Then
$$\sum_{k=0}^{p-1}\f{\b{2k}k\b{3k}k\b{6k}{3k}}{(-960)^{3k}}
\e\cases 0\mod {p^2}&\t{if $\sls p{43}=-1$,}
\\\sls p{15}x^2\mod p&\t{if $\sls p{43}=1$ and so $4p=x^2+43y^2$.}
\endcases$$
\endpro
Proof. Taking $m=-80\cdot 43$ and $n=42\cdot 43^2$ in Theorem 4.2
and then applying (2.9) we deduce the result.

\pro{Theorem 4.6 ([Su4, Conjecture 2.9])} Let $p$ be a prime such
that $p\not=2,3,5,11,67$. Then
$$\sum_{k=0}^{p-1}\f{\b{2k}k\b{3k}k\b{6k}{3k}}{(-5280)^{3k}}
\e\cases 0\mod {p^2}&\t{if $\sls p{67}=-1$,}
\\\sls {-330}px^2\mod p&\t{if $\sls p{67}=1$ and so $4p=x^2+67y^2$.}
\endcases$$
\endpro
Proof. Taking  $m=-440\cdot 67$ and $n=434\cdot 67^2$ in Theorem 4.2
and then applying (2.9) we deduce the result.

\pro{Theorem 4.7 ([Su4, Conjecture 2.10])} Let $p$ be a prime with
$p\not=2,3,5,23,29,163$. Then
$$\sum_{k=0}^{p-1}\f{\b{2k}k\b{3k}k\b{6k}{3k}}{(-640320)^{3k}}
\e\cases 0\mod {p^2}&\t{if $\sls p{163}=-1$,}
\\\sls {-10005}px^2\mod p&\t{if $\sls p{163}=1$ and so $4p=x^2+163y^2$.}
\endcases$$
\endpro
Proof. Taking  $m=-80\cdot 23\cdot 29\cdot 163$ and $n=14\cdot
11\cdot 19\cdot 127\cdot 163^2$ in Theorem 4.2 and then applying
(2.9) we deduce the result.

\pro{Theorem 4.8 ([S4, Conjecture 2.8])} Let $p>7$ be a prime. Then
$$\sum_{k=0}^{p-1}\f{\b{2k}k\b{3k}k\b{6k}{3k}}{(-15)^{3k}}
\e\cases 0\mod {p^2}&\t{if $p\e 3,5,6\mod 7$,}
\\\sls p{15}4C^2\mod p&\t{if $p=C^2+7D^2\e 1,2,4\mod 7$}.
\endcases$$
\endpro
Proof. Taking  $m=-35$ and $n=98$ in Theorem 4.2 and then applying
(2.5) we deduce the result.

\pro{Theorem 4.9 ([S4, Conjecture 2.9])} Let $p>7$ be a prime and
$p\not=17$. Then
$$\sum_{k=0}^{p-1}\f{\b{2k}k\b{3k}k\b{6k}{3k}}{255^{3k}} \e\cases
0\mod {p^2}&\t{if $p\e 3,5,6\mod 7$,}
\\\sls p{255}4C^2\mod p&\t{if $p=C^2+7D^2\e 1,2,4\mod 7$.}
\endcases$$
\endpro
Proof. Taking  $m=-595$ and $n=5586$ in Theorem 4.2 and then
applying (2.7) we deduce the result.

\pro{Theorem 4.10 ([S4, Conjecture 2.4])} Let $p$ be a prime such
that $p\not=2,3,11$. Then
$$\sum_{k=0}^{p-1}\f{\b{2k}k\b{3k}k\b{6k}{3k}}{66^{3k}}
\e\cases 0\mod {p^2}&\t{if $p\e 3\mod 4$,}
\\\sls p{33}4a^2\mod p&\t{if $p=a^2+b^2\e 1\mod 4$ and $2\nmid a$.}
\endcases$$
\endpro
Proof. Taking  $m=-11$ and $n=14$ in Theorem 4.2 and then applying
(2.2) we deduce the result.

\pro{Theorem 4.11 ([S4, Conjecture 2.5])} Let $p>5$ be a prime. Then
$$\sum_{k=0}^{p-1}\f{\b{2k}k\b{3k}k\b{6k}{3k}}{20^{3k}}
\e\cases 0\mod {p^2}&\t{if $p\e 5,7\mod 8$,}
\\\sls {-5}p4c^2\mod p&\t{if $p=c^2+2d^2\e 1,3\mod 8$}.
\endcases$$
\endpro
Proof. Taking  $m=-30$ and $n=56$ in Theorem 4.2 and then applying
the result in the proof of Corollary 2.2 we deduce the result.

\pro{Theorem 4.12 ([S4, Conjecture 2.6])} Let $p>5$ be a prime. Then
$$\sum_{k=0}^{p-1}\f{\b{2k}k\b{3k}k\b{6k}{3k}}{54000^k}
\e\cases 0\mod {p^2}&\t{if $p\e 2\mod 3$,}
\\\sls p54A^2\mod p&\t{if $p=A^2+3B^2\e 1\mod 3$.}
\endcases$$
\endpro
Proof. Taking  $m=-15$ and $n=22$ in Theorem 4.2 and then applying
(2.3) we deduce the result.
 \pro{Theorem 4.13 ([S4, Conjecture 2.7])} Let $p>5$ be a prime.
Then
$$\sum_{k=0}^{p-1}\f{\b{2k}k\b{3k}k\b{6k}{3k}}{(-12288000)^k}
\e\cases 0\mod{p^2}&\t{if $3\mid p-2$,}\\\sls {10}pL^2\mod p &\t{if
$3\mid p-1$ and so $4p=L^2+27M^2$.}
\endcases$$
\endpro
Proof. Taking  $m=-120$ and $n=506$ in Theorem 4.2 and then applying
(2.4) we deduce the result.
\newline{\bf Remark 4.1} From [O] we
know that the only $j$-invariants of elliptic curves over rational
field $\Bbb Q$ with complex multiplication are given by
$$0,12^3,-15^3,20^3,-32^3,2\cdot 30^3,66^3,-96^3,-3\cdot 160^3,
255^3,-960^3,-5280^3,-640320^3,$$ coinciding with the values of $m$
 in (3.3) and Theorems 4.3-4.13.

 \subheading{5. Some conjectures on
supercongruences}
 \pro{Conjecture 5.1} Let $p>5$ be a prime. Then
$$\align &\sum_{k=0}^{p-1}\f{63k+8}{(-15)^{3k}}
\b{2k}k\b{3k}k\b{6k}{3k}\e 8p\Ls {-15}p\mod{p^2},
\\&\sum_{k=0}^{p-1}\f{133k+8}{255^{3k}}
\b{2k}k\b{3k}k\b{6k}{3k}\e 8p\Ls {-255}p\mod{p^2}\qtq{for}p\not=17,
\\&\sum_{k=0}^{p-1}\f{28k+3}{20^{3k}}
\b{2k}k\b{3k}k\b{6k}{3k}\e 3p\Ls{-5}p\mod{p^2},
\\&\sum_{k=0}^{p-1}\f{63k+5}{66^{3k}}
\b{2k}k\b{3k}k\b{6k}{3k}\e 5p\Ls{-33}p\mod{p^2}\qtq{for}p\not=11,
\\&\sum_{k=0}^{p-1}\f{11k+1}{54000^k}
\b{2k}k\b{3k}k\b{6k}{3k}\e p\Ls {-15}p\mod{p^2},
\\&\sum_{k=0}^{p-1}\f{506k+31}{(-12288000)^k}
\b{2k}k\b{3k}k\b{6k}{3k}\e 31p\Ls{-30}p\mod{p^2}.
\endalign$$
\endpro
\par Conjecture 5.1 is similar to some conjectures in [Su1].
\pro{Conjecture 5.2} Let $p>3$ be a prime. Then
$$\align&\sum_{n=0}^{p-1}\f{9n+4}{5^n}\b{2n}n\sum_{k=0}^n\b nk^3
\e 4p\Ls p5\mod{p^2}\qtq{for}p>5,
\\&\sum_{n=0}^{p-1}\f{5n+2}{16^n}\b{2n}n\sum_{k=0}^n\b nk^3
\e 2p\mod{p^2},
\\&\sum_{n=0}^{p-1}\f{9n+2}{50^n}\b{2n}n\sum_{k=0}^n\b nk^3
\e 2p\Ls {-1}p\mod{p^2}\qtq{for}p\not=5,
\\&\sum_{n=0}^{p-1}\f{5n+1}{96^n}\b{2n}n\sum_{k=0}^n\b nk^3
\e p\Ls {-2}p\mod{p^2},
\\&\sum_{n=0}^{p-1}\f{6n+1}{320^n}\b{2n}n\sum_{k=0}^n\b nk^3
\e p\Ls p{15}\mod{p^2}\qtq{for}p\not=5,
\\&\sum_{n=0}^{p-1}\f{90n+13}{896^n}\b{2n}n\sum_{k=0}^n\b nk^3
\e 13p\Ls p{7}\mod{p^2}\qtq{for}p\not=7,
\\&\sum_{n=0}^{p-1}\f{102n+11}{10400^n}\b{2n}n\sum_{k=0}^n\b nk^3
\e 11p\Ls p{39}\mod{p^2}\qtq{for}p\not=5,13.
\endalign$$
\endpro

\pro{Conjecture 5.3} Let $p>3$ be a prime. Then
$$\align
\\&\sum_{n=0}^{p-1}\f{3n+1}{(-16)^n}\b{2n}n\sum_{k=0}^n\b nk^3
\e p\Ls {-1}p\mod{p^2},
\\&\sum_{n=0}^{p-1}\f{15n+4}{(-49)^n}\b{2n}n\sum_{k=0}^n\b nk^3
\e 4p\Ls p3\mod{p^2}\qtq{for}p\not=7,
\\&\sum_{n=0}^{p-1}\f{9n+2}{(-112)^n}\b{2n}n\sum_{k=0}^n\b nk^3
\e 2p\Ls p7\mod{p^2}\qtq{for}p\not=7,
\\&\sum_{n=0}^{p-1}\f{99n+17}{(-400)^n}\b{2n}n\sum_{k=0}^n\b nk^3
\e 17p\Ls{-1}p\mod{p^2},
\\&\sum_{n=0}^{p-1}\f{855n+109}{(-2704)^n}\b{2n}n\sum_{k=0}^n\b nk^3
\e 109p\Ls{-1}p\mod{p^2}\qtq{for}p\not=13,
\\&\sum_{n=0}^{p-1}\f{585n+58}{(-24304)^n}\b{2n}n\sum_{k=0}^n\b nk^3
\e 58p\Ls{-31}p\mod{p^2}\qtq{for}p\not=7,31.
\endalign$$
\endpro

\par For an integer $m$ and odd prime $p$ with $p\nmid m$ let

$$Z_p(m)=\sum_{n=0}^{p-1}\f{\b{2n}n}{m^n}\sum_{k=0}^n\b nk^3.$$
Then we have the following conjectures concerning $Z_p(m)$ modulo
$p^2$.
 \pro{Conjecture 5.4} Let $p$ be an odd prime. Then
$$Z_p(-16)\e \cases 4x^2-2p\mod {p^2}&\t{if $p=x^2+y^2\e 1\mod{12}$ with $6\mid y$,}
\\2p-4x^2\mod{p^2}&\t{if $p=x^2+y^2\e 1\mod{12}$ with $6\mid x-3$,}
\\4\sls {xy}3xy\mod{p^2}&\t{if $p=x^2+y^2\e 5\mod{12}$,}
\\0\mod{p^2}&\t{if $p\e 3\mod 4$.}\endcases
$$\endpro
 \pro{Conjecture 5.5} Let $p$ be an odd prime. Then
$$Z_p(96)\e \cases \sls p3(4x^2-2p)\mod{p^2}&\t{if $p=x^2+2y^2\e
1,3\mod 8$,}\\0\mod{p^2}&\t{if $p\e 5,7\mod 8$.}
\endcases$$
\endpro
\pro{Conjecture 5.6} Let $p>5$ be a prime. Then
$$Z_p(-4)\e Z_p(50)\e \cases 4x^2-2p\mod{p^2}&\t{if $p=x^2+3y^2\e
1\mod 3$,}\\0\mod{p^2}&\t{if $p\e 2\mod 3$.}
\endcases$$
\endpro
\pro{Conjecture 5.7} Let $p>5$ be a prime. Then
$$Z_p(16)\e \cases 4x^2-2p\mod{p^2}&\t{if $p\e
1,9\mod{20}$ and so $p=x^2+5y^2$,}\\2x^2-2p\mod{p^2}&\t{if $p\e
3,7\mod{20}$ and so $2p=x^2+5y^2$,}
\\0\mod{p^2}&\t{if $p\e 11,13,17,19\mod{20}$.}\endcases$$\endpro

\pro{Conjecture 5.8} Let $p>3$ be a prime. Then
$$Z_p(32)\e \cases 4x^2-2p\mod{p^2}&\t{if $p\e
1,7\mod{24}$ and so $p=x^2+6y^2$,}\\8x^2-2p\mod{p^2}&\t{if $p\e
5,11\mod{24}$ and so $p=2x^2+3y^2$,}\\0\mod{p^2}&\t{if $p\e
13,17,19,23\mod{24}$.}\endcases$$\endpro

\pro{Conjecture 5.9} Let $p>7$ be a prime. Then
$$Z_p(5)\e Z_p(-49)\e \cases 4x^2-2p\mod{p^2}&\t{if $p=x^2+15y^2\e
1,19\mod{30}$,}\\2p-12x^2\mod{p^2}&\t{if $p=3x^2+5y^2\e
17,23\mod{30}$,}
\\0\mod{p^2}&\t{if $p\e 7,11,13,29\mod{30}$.}\endcases$$\endpro

\pro{Conjecture 5.10} Let $b\in\{7,11,19,31,59\}$ and let
$f(b)=-112,-400,-2704$, $-24304$, $-1123600$ according as $b=
7,11,19,31,59$. If $p$ is a prime with $p\not= 2,3,b$ and $p\nmid
f(b)$, then
$$Z_p(f(b))\e \cases 4x^2-2p\mod{p^2}&\t{if $p=x^2+3by^2$,}
\\2p-12x^2\mod{p^2}&\t{if $p=3x^2+by^2$,}
\\2x^2-2p\mod{p^2}&\t{if $2p=x^2+3by^2$,}
\\2p-6x^2\mod{p^2}&\t{if $2p=3x^2+by^2$,}
\\0\mod{p^2}&\t{if $\sls{-3b}p=-1$.}
\endcases$$\endpro

\pro{Conjecture 5.11} Let $b\in\{5,7,13,17\}$ and $f(b)= 320,
896,10400,39200$ according as $b=5,7,13,17$. If $p$ is a prime with
$p\not= 2,3,b$ and $p\nmid f(b)$, then
$$Z_p(f(b))\e \cases 4x^2-2p\mod{p^2}&\t{if $p=x^2+6by^2$,}
\\8x^2-2p\mod{p^2}&\t{if $p=2x^2+3by^2$,}
\\2p-12x^2\mod{p^2}&\t{if $p=3x^2+2by^2$,}
\\2p-24x^2\mod{p^2}&\t{if $p=6x^2+by^2$,}
\\0\mod{p^2}&\t{if $\sls{-6b}p=-1$.}
\endcases$$\endpro

\enddocument
\bye